\documentclass{article}
\usepackage{amsfonts}

\newtheorem{theorem}{Theorem}

\newtheorem{corollary}[theorem]{Corollary}

\newtheorem{proposition}[theorem]{Proposition}
\newtheorem{remark}[theorem]{Remark}

\newenvironment{proof}[1][Proof]{\noindent\textbf{#1.} }{\ \rule{0.5em}{0.5em}}
\input{tcilatex}

\begin{document}

\title{Confidence Sets Based on Sparse Estimators Are Necessarily Large}
\author{Benedikt M. P\"{o}tscher\thanks{%
Department of Statistics, University of Vienna, Universit\"{a}tsstrasse 5,
A-1010 Vienna. Phone: +431 427738640. E-mail: benedikt.poetscher@univie.ac.at%
} \\
Department of Statistics, University of Vienna}
\date{Preliminary version: April 2007\\
First version: August 2007\\
Revised version: April 2009}
\maketitle

\begin{abstract}
Confidence sets based on sparse estimators are shown to be large compared to
more standard confidence sets, demonstrating that sparsity of an estimator
comes at a substantial price in terms of the quality of the estimator. The
results are set in a general parametric or semiparametric framework.


62J07


post-model-selection estimator, penalized maximum likelihood,

confidence set, coverage probability
\end{abstract}

\section{Introduction}

Sparse estimators have received increased attention in the statistics
literature in recent years. An estimator for a parameter vector is called
sparse if it estimates the zero components of the true parameter vector by
zero with probability approaching one as sample size increases without
bound. Examples of sparse estimator are (i) post-model-selection estimators
following a consistent model selection procedure, (ii) thresholding
estimators with a suitable choice of the thresholds, and (iii) many
penalized maximum likelihood estimators (e.g., SCAD, LASSO, and variants
thereof)\ when the regularization parameter is chosen in a suitable way.
Many (but not all) of these sparse estimators also have the property that
the asymptotic distribution of the estimator coincides with the asymptotic
distribution of the (infeasible) estimator that uses the zero restrictions
in the true parameter; see, e.g., P\"{o}tscher (1991, Lemma 1), Fan and Li
(2001). This property has -- in the context of SCAD estimation -- been
dubbed the \textquotedblleft oracle\textquotedblright\ property by Fan and
Li (2001) and has received considerable attention in the literature,
witnessed by a series of papers establishing the \textquotedblleft
oracle\textquotedblright\ property for a variety of estimators (e.g., Bunea
(2004), Bunea and McKeague (2005), Fan and Li (2002, 2004), Zou (2006), Li
and Liang (2007), Wang and Leng (2007), Wang, G. Li, and Tsai (2007), Wang,
R. Li, and Tsai (2007), Zhang and Lu (2007), Zou and Yuan (2008)).

The sparsity property and the closely related \textquotedblleft
oracle\textquotedblright\ property seem to intimate that an estimator
enjoying these properties is superior to classical estimators like the
maximum likelihood estimator (not possessing the \textquotedblleft
oracle\textquotedblright\ property). We show, however, that the sparsity
property of an estimator does not translate into good properties of
confidence sets based on this estimator. Rather we show in Section 2 that
any confidence set based on a sparse estimator is necessarily large relative
to more standard confidence sets, e.g., obtained from the maximum likelihood
estimator, that have the same guaranteed coverage probability. Hence, there
is a substantial price to be paid for sparsity, which is not revealed by the
pointwise asymptotic analysis underlying the \textquotedblleft
oracle\textquotedblright\ property. Special cases of the general results
provided in Section 2 have been observed in the literature: It has been
noted that the \textquotedblleft naive\textquotedblright\ confidence
interval centered at Hodges' estimator has infimal coverage probability that
converges to zero as sample size goes to infinity, see Kale (1985), Beran
(1992), and Kabaila (1995). [By the \textquotedblleft
naive\textquotedblright\ confidence interval we mean the interval one would
construct in the usual way from the pointwise asymptotic distribution of
Hodges' estimator.] Similar results for \textquotedblleft
naive\textquotedblright\ confidence intervals centered at
post-model-selection estimators that are derived from certain consistent
model selection procedures can be found in Kabaila (1995) and Leeb and P\"{o}%
tscher (2005). We note that these \textquotedblleft naive\textquotedblright\
confidence intervals have coverage probabilities that converge to the
nominal level \emph{pointwise }in the parameter space, but these confidence
intervals are -- in view of the results just mentioned -- not
\textquotedblleft honest\textquotedblright\ in the sense that the infimum
over the parameter space of the coverage probabilities converges to a level
that is \emph{below} the nominal level. Properties of confidence sets based
on not necessarily sparsely tuned post-model-selection estimators are
discussed in Kabaila (1995, 1998), P\"{o}tscher (1995), Leeb and P\"{o}%
tscher (2005), Kabaila and Leeb (2006).

The results discussed in the preceding paragraph show, in particular, that
the \textquotedblleft oracle\textquotedblright\ property is problematic as
it gives a much too optimistic impression of the actual properties of an
estimator. This problematic nature of the \textquotedblleft
oracle\textquotedblright\ property is also discussed in Leeb and P\"{o}%
tscher (2008) from a risk point of view; cf. also Yang (2005). The
problematic nature of the \textquotedblleft oracle\textquotedblright\
property is connected to the fact that the finite-sample distributions of
these estimators converge to their limits pointwise in the parameter space
but not uniformly. Hence, the limits often do not reveal the actual
properties of the finite-sample distributions. An asymptotic analysis using
a "moving parameter" asymptotics is possible and captures much of the actual
behavior of the estimators, see Leeb and P\"{o}tscher (2005), P\"{o}tscher
and Leeb (2007), and P\"{o}tscher and Schneider (2009). These results lead
to a view of these estimators that is less favorable then what is suggested
by the \textquotedblleft oracle\textquotedblright\ property.

The remainder of the paper is organized as follows: In Section 2 we provide
the main results showing that confidence sets based on sparse estimators are
necessarily large. These results are extended to \textquotedblleft
partially\textquotedblright\ sparse estimators in Section 2.1. In Section 3
we consider a thresholding estimator as a simple example of a sparse
estimator, construct a confidence set based on this estimator, and discuss
its properties.

\section{On the size of confidence sets based on sparse estimators}

Suppose we are given a sequence of statistical experiments 
\begin{equation}
\left\{ P_{n,\theta }:\theta \in \mathbb{R}^{k}\right\} \qquad n=1,2,\ldots
\label{0}
\end{equation}%
where the probability measures $P_{n,\theta }$ live on suitable measure
spaces $(\mathcal{X}_{n},\mathfrak{X}_{n})$. [Often $P_{n,\theta }$ will
arise as the distribution of a random vector $(y_{1}^{\prime },\ldots
,y_{n}^{\prime })^{\prime }$ where $y_{i}$ takes values in a Euclidean
space. In this case $\mathcal{X}_{n}$ will be an $n$-fold product of that
Euclidean space and $\mathfrak{X}_{n}$ will be the associated Borel $\sigma $%
-field; also $n$ will then denote sample size.] We assume further that for
every $\gamma \in \mathbb{R}^{k}$ the sequence of probability measures 
\[
\left\{ P_{n,\gamma /\sqrt{n}}:n=1,2,\ldots \right\} 
\]%
is contiguous w.r.t. the sequence 
\[
\left\{ P_{n,0}:n=1,2,\ldots \right\} . 
\]%
This is a quite weak assumption satisfied by many statistical experiments
(including experiments with dependent data); for example, it is certainly
satisfied whenever the experiment is locally asymptotically normal. The
above assumption that the parameter space is $\mathbb{R}^{k}$ is made only
for simplicity of presentation and is by no means essential, see Remark \ref%
{parameterspace}.

Let $\hat{\theta}_{n}$ denote a sequence of estimators, i.e., $\hat{\theta}%
_{n}$ is a measurable function on $\mathcal{X}_{n}$ taking values in $%
\mathbb{R}^{k}$. We say that the estimator $\hat{\theta}_{n}$\ (more
precisely, the sequence of estimators) is \textit{sparse }if for every $%
\theta \in \mathbb{R}^{k}$ and $i=1,\ldots ,k$%
\begin{equation}
\lim_{n\rightarrow \infty }P_{n,\theta }\left( \hat{\theta}_{n,i}=0\right) =1%
\text{ \ \ \ holds whenever \ \ \ }\theta _{i}=0\text{.}  \label{1}
\end{equation}%
Here $\hat{\theta}_{n,i}$ and $\theta _{i}$ denote the $i$-th component of $%
\hat{\theta}_{n}$ and of $\theta $, respectively. That is, the estimator is
guaranteed to find the zero components of $\theta $ with probability
approaching one as $n\rightarrow \infty $. [The focus on zero-values in the
coordinates of $\theta $ is of course arbitrary. Furthermore, note that
Condition (\ref{1}) is of course satisfied for nonsensical estimators like $%
\hat{\theta}_{n}\equiv 0$. The sparse estimators mentioned in Section 1 and
Remark \ref{example} below, however, are more sensible as they are typically
also consistent for $\theta $.]

\begin{remark}
\label{example}Typical examples of sparse estimators are as follows:
consider a linear regression model $Y=X\theta +u$ under standard assumptions
(for simplicity assume $u\sim N(0,\sigma ^{2}I_{n})$ with $\sigma >0$ known
and $X$ nonstochastic with $X^{\prime }X/n\rightarrow Q$, a positive
definite matrix). Suppose a subset of the regressors contained in $X$ is
selected first by an application of a consistent all-subset model selection
procedure (such as, e.g., Schwarz' minimum BIC-method) and then the least
squares estimator based on the selected model is reported, with the
coefficients of the excluded regressor variables being estimated as zero.
The resulting estimator for $\theta $ is a so-called post-model-selection
estimator and clearly has the sparsity property. Another estimator
possessing the sparsity property can be obtained via hard-thresholding as
follows: compute the least squares estimator from the full model $Y=X\theta
+u$ and replace those components of the least squares estimator by zero
which have a $t$-statistic that is less than a threshold $\eta _{n}$ in
absolute value. The resulting estimator has the sparsity property if $\eta
_{n}\rightarrow 0$ and $n^{1/2}\eta _{n}\rightarrow \infty $ holds for $%
n\rightarrow \infty $. As mentioned in the introduction, also a large class
of penalized least squares estimators has the sparsity property, see the
references given there.
\end{remark}

Returning to the general discussion, we are interested in confidence sets
for $\theta $ based on $\hat{\theta}_{n}$. Let $C_{n}$ be a random set in $%
\mathbb{R}^{k}$ in the sense that $C_{n}=C_{n}(\omega )$ is a subset of $%
\mathbb{R}^{k}$ for every $\omega \in \mathcal{X}_{n}$ with the property
that for every $\theta \in \mathbb{R}^{k}$%
\[
\left\{ \omega \in \mathcal{X}_{n}:\theta \in C_{n}(\omega )\right\} 
\]%
is measurable, i.e., belongs to $\mathfrak{X}_{n}$. We say that the random
set $C_{n}$ is \textit{based on }the estimator $\hat{\theta}_{n}$ if $C_{n}$
satisfies%
\begin{equation}
P_{n,\theta }\left( \hat{\theta}_{n}\in C_{n}\right) =1  \label{2}
\end{equation}%
for every $\theta \in \mathbb{R}^{k}$. [If the set inside of the probability
in (\ref{2}) is not measurable, the probability is to be replaced by inner
probability.]\ For example, if $C_{n}$ is a $k$-dimensional interval (box)
of the form 
\begin{equation}
\left[ \hat{\theta}_{n}-a_{n},\hat{\theta}_{n}+b_{n}\right]  \label{box}
\end{equation}%
where $a_{n}$ and $b_{n}$ are random vectors in $\mathbb{R}^{k}$ with only
nonnegative coordinates, then condition (\ref{2}) is trivially satisfied.
Here we use the notation $[c,d]=[c_{1},d_{1}]\times \cdots \times \lbrack
c_{k},d_{k}]$ for vectors $c=(c_{1},\ldots ,c_{k})^{\prime }$ and $%
d=(d_{1},\ldots ,d_{k})^{\prime }$. We also use the following notation: For
a subset $A$ of $\mathbb{R}^{k}$, let 
\[
\limfunc{diam}(A)=\sup \{\left\Vert x-y\right\Vert :x\in A,y\in A\} 
\]%
denote the diameter of $A$ (measured w.r.t. the usual Euclidean norm $%
\left\Vert \cdot \right\Vert $); furthermore, if $e$ is an arbitrary element
of $\mathbb{R}^{k}$ of length $1$, and $a\in A$ let 
\[
\limfunc{ext}(A,a,e)=\sup \{\lambda \geq 0:\lambda e+a\in A\}. 
\]%
That is, $\limfunc{ext}(A,a,e)$ measures how far the set $A$ extends from
the point $a$ into the direction given by $e$. [Observe that without further
conditions (such as, e.g., convexity of $A$) not all points of the form $%
\lambda e+a$ with $\lambda <\limfunc{ext}(A,a,e)$ need to belong to $A$.]

The following result shows that confidence sets based on a sparse estimator
are necessarily large.

\begin{theorem}
\label{Th1}Suppose the statistical experiment given in (\ref{0}) satisfies
the above contiguity assumption. Let $\hat{\theta}_{n}$ be a sparse
estimator sequence and let $C_{n}$ be a sequence of random sets \textit{%
based on }the estimator $\hat{\theta}_{n}$ in the sense of (\ref{2}). Assume
that $C_{n}$ is a confidence set for $\theta $ with asymptotic infimal
coverage probability $\delta $, i.e.,%
\[
\delta =\liminf_{n\rightarrow \infty }\inf_{\theta \in \mathbb{R}%
^{k}}P_{n,\theta }\left( \theta \in C_{n}\right) . 
\]%
Then for every $t\geq 0$ and every $e\in \mathbb{R}^{k}$ of length $1$ we
have%
\begin{equation}
\liminf_{n\rightarrow \infty }\sup_{\theta \in \mathbb{R}^{k}}P_{n,\theta
}\left( \sqrt{n}\limfunc{ext}(C_{n},\hat{\theta}_{n},e)\geq t\right) \geq
\delta .  \label{ext}
\end{equation}%
In particular, we have for every $t\geq 0$%
\begin{equation}
\liminf_{n\rightarrow \infty }\sup_{\theta \in \mathbb{R}^{k}}P_{n,\theta
}\left( \sqrt{n}\limfunc{diam}(C_{n})\geq t\right) \geq \delta .
\label{diam}
\end{equation}%
[If the set inside of the probability in (\ref{ext}) or (\ref{diam}) is not
measurable, the probability is to be replaced by inner probability.]
\end{theorem}

\begin{proof}
Since obviously $\limfunc{diam}(C_{n})\geq \limfunc{ext}(C_{n},\hat{\theta}%
_{n},e)$ holds with $P_{n,\theta }$-probability $1$ for all $\theta $ in
view of (\ref{2}), it suffices to prove (\ref{ext}). Now, for every sequence 
$\theta _{n}\in \mathbb{R}^{k}$ we have in view of (\ref{2}) 
\begin{eqnarray}
\delta &=&\liminf_{n\rightarrow \infty }\inf_{\theta \in \mathbb{R}%
^{k}}P_{n,\theta }\left( \theta \in C_{n}\right) \leq \liminf_{n\rightarrow
\infty }P_{n,\theta _{n}}\left( \theta _{n}\in C_{n}\right)  \label{3} \\
&=&\liminf_{n\rightarrow \infty }\left\{ P_{n,\theta _{n}}\left( \theta
_{n}\in C_{n},\hat{\theta}_{n}\in C_{n},\hat{\theta}_{n}=0\right) \right. 
\nonumber \\
&&\left. +P_{n,\theta _{n}}\left( \theta _{n}\in C_{n},\hat{\theta}_{n}\neq
0\right) \right\} .  \nonumber
\end{eqnarray}%
Sparsity implies%
\[
\lim_{n\rightarrow \infty }P_{n,0}\left( \hat{\theta}_{n}\neq 0\right) =0, 
\]%
and hence for $\theta _{n}=\gamma /\sqrt{n}$ the contiguity assumption
implies%
\[
\limsup_{n\rightarrow \infty }P_{n,\theta _{n}}\left( \theta _{n}\in C_{n},%
\hat{\theta}_{n}\neq 0\right) \leq \lim_{n\rightarrow \infty }P_{n,\theta
_{n}}\left( \hat{\theta}_{n}\neq 0\right) =0. 
\]%
Consequently, we obtain from (\ref{3}) for $\theta _{n}=\gamma /\sqrt{n}$
with $\gamma \neq 0$%
\begin{eqnarray}
\delta &\leq &\liminf_{n\rightarrow \infty }P_{n,\theta _{n}}\left( \theta
_{n}\in C_{n},\hat{\theta}_{n}\in C_{n},\hat{\theta}_{n}=0\right)  \nonumber
\\
&\leq &\liminf_{n\rightarrow \infty }P_{n,\theta _{n}}\left( \sqrt{n}%
\limfunc{ext}(C_{n},\hat{\theta}_{n},\gamma /\left\Vert \gamma \right\Vert
)\geq \left\Vert \gamma \right\Vert \right)  \label{4}
\end{eqnarray}%
because of the obvious inclusion%
\[
\left\{ \theta _{n}\in C_{n},\hat{\theta}_{n}\in C_{n},\hat{\theta}%
_{n}=0\right\} \subseteq \left\{ \limfunc{ext}(C_{n},\hat{\theta}_{n},\theta
_{n}/\left\Vert \theta _{n}\right\Vert )\geq \left\Vert \theta
_{n}\right\Vert \right\} . 
\]%
Since $\gamma $ was arbitrary, the result (\ref{ext})\ follows from (\ref{4}%
) upon identifying $t$ and $\left\Vert \gamma \right\Vert $.
\end{proof}

\begin{corollary}
Suppose the assumptions of Theorem \ref{Th1} are satisfied and $C_{n}$ is a
confidence `interval' of the form (\ref{box}). Then for every $i=1,\ldots ,k$
and every $t\geq 0$ 
\[
\liminf_{n\rightarrow \infty }\sup_{\theta \in \mathbb{R}^{k}}P_{n,\theta
}\left( \sqrt{n}a_{n,i}\geq t\right) \geq \delta 
\]%
and%
\[
\liminf_{n\rightarrow \infty }\sup_{\theta \in \mathbb{R}^{k}}P_{n,\theta
}\left( \sqrt{n}b_{n,i}\geq t\right) \geq \delta 
\]%
hold, where $a_{n,i}$ and $b_{n,i}$ denote the $i$-th coordinate of $a_{n}$
and $b_{n}$, respectively. In particular, if $a_{n}$ and $b_{n}$ are
nonrandom, 
\[
\liminf_{n\rightarrow \infty }\sqrt{n}a_{n,i}=\liminf_{n\rightarrow \infty }%
\sqrt{n}b_{n,i}=\infty 
\]%
holds for every $i=1,\ldots ,k$, provided that $\delta >0$.
\end{corollary}

\begin{proof}
Follows immediately from the previous theorem upon observing that (\ref{box}%
) implies $\limfunc{ext}(C_{n},\hat{\theta}_{n},-e_{i})=a_{n,i}$ and $%
\limfunc{ext}(C_{n},\hat{\theta}_{n},e_{i})=b_{n,i}$ where $e_{i}$ denotes
the $i$-th standard basis vector.
\end{proof}

It is instructive to compare with standard confidence sets. For example, in
a normal linear regression model $\sqrt{n}$ times the diameter of the
standard confidence ellipsoid is stochastically bounded uniformly in $\theta 
$. In contrast, Theorem \ref{Th1} tells us that any confidence set $C_{n}$
based on sparse estimators with $\sqrt{n}\limfunc{diam}(C_{n})$ being
stochastically bounded uniformly in $\theta $ necessarily has infimal
coverage probability equal to zero.

\begin{remark}
\label{nuis}(Nuisance parameters) Suppose that the sequence of statistical
experiments is of the form $\left\{ P_{n,\theta ,\tau }:\theta \in \mathbb{R}%
^{k},\tau \in T\right\} $ where $\theta $ is the parameter of interest and $%
\tau $ is now a (possibly infinite dimensional) nuisance parameter. Theorem %
\ref{Th1} can then clearly be applied to the parametric subfamilies $\left\{
P_{n,\theta ,\tau }:\theta \in \mathbb{R}^{k}\right\} $ for $\tau \in T$
(provided the conditions of the theorem are satisfied). In particular, the
following is then an immediate consequence: suppose that the contiguity
condition and sparsity condition are satisfied for every $\tau \in T$.
Suppose further that we are again interested in confidence sets for $\theta $
based on $\hat{\theta}_{n}$ (in the sense that $P_{n,\theta ,\tau }\left( 
\hat{\theta}_{n}\in C_{n}\right) =1$ for all $\theta \in \mathbb{R}^{k},\tau
\in T$) that have asymptotic infimal (over $\theta $ and $\tau $) coverage
probability $\delta $. Then results analogous to (\ref{ext}) and (\ref{diam}%
), but with the supremum extending now over $\mathbb{R}^{k}\times T$, hold.
\end{remark}

\begin{remark}
\label{funct}(Confidence sets for linear functions of $\theta $)\ Suppose
that a statistical experiment $\left\{ P_{n,\theta }:\theta \in \mathbb{R}%
^{k}\right\} $ satisfying the aforementioned contiguity property and a
sparse estimator $\hat{\theta}_{n}$ are given but that we are interested in
setting a confidence set for $\vartheta =A\theta $ that is based on $\hat{%
\vartheta}_{n}=A\hat{\theta}_{n}$, where $A$ is a given $q\times k$ matrix.
Without loss of generality assume that $A$ has full row rank. [In
particular, this covers the case where we have a sparse estimator for $%
\theta $, but are interested in confidence sets for a subvector only.]
Suppose $C_{n}$ is a confidence set for $\vartheta $ that is based on $\hat{%
\vartheta}_{n}$ (in the sense that $P_{n,\theta }\left( \hat{\vartheta}%
_{n}\in C_{n}\right) =1$ for all $\theta \in \mathbb{R}^{k}$) and that has
asymptotic infimal coverage probability $\delta $. Then essentially the same
proof as for Theorem \ref{Th1} shows that for every $t\geq 0$ and every $%
e\in \mathbb{R}^{q}$ of length $1$ we have%
\begin{equation}
\liminf_{n\rightarrow \infty }\sup_{\theta \in \mathbb{R}^{k}}P_{n,\theta
}\left( \sqrt{n}\limfunc{ext}(C_{n},\hat{\vartheta}_{n},e)\geq t\right) \geq
\delta  \label{lin_funct}
\end{equation}%
and consequently also the analogue of (\ref{diam}) holds.
\end{remark}

\begin{remark}
\label{rate}The contiguity assumption together with the sparsity of the
estimator was used in the proof of Theorem \ref{Th1} to imply $%
\lim_{n\rightarrow \infty }P_{n,\theta _{n}}\left( \hat{\theta}_{n}\neq
0\right) =0$ for all sequences of the form $\theta _{n}=\gamma /\sqrt{n}$, $%
\gamma \in \mathbb{R}^{k}$. For some important classes of sparse estimators
this relation can even be established for all sequences of the form $\theta
_{n}=\gamma /v_{n}$, $\gamma \in \mathbb{R}^{k}$, where $v_{n}$ are certain
sequences that diverge to infinity, but at a rate slower than $\sqrt{n}$
(cf. Leeb and P\"{o}tscher (2005), P\"{o}tscher and Leeb (2007, Proposition
1), P\"{o}tscher and Schneider (2009, Proposition 1)). Inspection of the
proof of Theorem \ref{Th1} shows that then a stronger result follows, namely
that (\ref{ext}) and (\ref{diam}) hold even with $\sqrt{n}$ replaced by $%
v_{n}$. This shows that in such a case confidence sets based on sparse
estimators are even larger than what is predicted by Theorem \ref{Th1}. This
simple extension immediately applies mutatis mutandis also to the other
results in the paper (with the exception of Theorem \ref{Th3}, an extension
of which would require a separate analysis). The example discussed in
Section 3 nicely illustrates the phenomenon just described.
\end{remark}

\begin{remark}
\label{parameterspace}The assumption that the parameter space indexing the
statistical experiment, say $\Theta $, is an entire Euclidean space is not
essential as can be seen from the proofs. The results equally well hold if,
e.g., $\Theta $ is a subset of Euclidean space that contains a ball with
center at zero (simply put $\theta _{n}=\gamma /\sqrt{n}$ if this belongs to 
$\Theta $, and set $\theta _{n}=0$ otherwise). In fact, $\Theta $ could even
be allowed to depend on $n$ and to \textquotedblleft
shrink\textquotedblright\ to zero at a rate slower than $n^{-1/2}$. [In that
sense the results are of a \textquotedblleft local\textquotedblright\ rather
than of a \textquotedblleft global\textquotedblright\ nature.]
\end{remark}

\begin{remark}
Suppose the contiguity assumption is satisfied and the estimator sequence $%
\hat{\theta}_{n}$ is sparse. Then the uniform convergence rate of $\hat{%
\theta}_{n}$ is necessarily slower than $n^{-1/2}$. In fact, more is true:
for every real number $M>0$ we have%
\begin{equation}
\liminf_{n\rightarrow \infty }\sup_{\theta \in \mathbb{R}^{k}}P_{n,\theta
}\left( n^{1/2}\left\Vert \hat{\theta}_{n}-\theta \right\Vert >M\right) =1.
\label{unif_rate}
\end{equation}%
To see this, set $\theta _{n}=\gamma /\sqrt{n}$ with $\left\Vert \gamma
\right\Vert >M$ and observe that the left-hand side in the above display is
not less than 
\[
\liminf_{n\rightarrow \infty }P_{n,\theta _{n}}\left( n^{1/2}\left\Vert \hat{%
\theta}_{n}-\theta _{n}\right\Vert >M\right) =\liminf_{n\rightarrow \infty
}P_{n,\theta _{n}}\left( n^{1/2}\left\Vert \theta _{n}\right\Vert >M,\hat{%
\theta}_{n}=0\right) =1,
\]%
the displayed equalities holding true in view of sparsity and contiguity.
[If $\lim_{n\rightarrow \infty }P_{n,\theta _{n}}\left( \hat{\theta}_{n}\neq
0\right) =0$ holds for all sequences of the form $\theta _{n}=\gamma /v_{n}$%
, $\gamma \in \mathbb{R}^{k}$, where $v_{n}>0$ is a given sequence (cf.
Remark \ref{rate}), then obviously (\ref{unif_rate}) holds with $n^{1/2}$
replaced by $v_{n}$. Furthermore, the results in this remark continue to
hold if the supremum over $\theta \in \mathbb{R}^{k}$ is replaced by a
supremum over a set $\Theta $ that contains a ball with center at zero.]
\end{remark}

\subsection{Confidence sets based on partially sparse estimators}

Suppose that in the framework of (\ref{0}) the parameter vector $\theta $ is
partitioned as $\theta =(\alpha ^{\prime },\beta ^{\prime })^{\prime }$
where $\alpha $ is $(k-k_{\beta })\times 1$ and $\beta $ is $k_{\beta
}\times 1$ ($0<k_{\beta }<k$). Furthermore, suppose that the estimator $\hat{%
\theta}_{n}=(\hat{\alpha}_{n}^{\prime },\hat{\beta}_{n}^{\prime })^{\prime }$
is `partially' sparse in the sense that it finds the zeros in $\beta $ with
probability approaching $1$ (but not necessarily the zeros in $\alpha $).
That is, for every $\theta \in \mathbb{R}^{k}$ and $i=1,\ldots ,k_{\beta }$%
\begin{equation}
\lim_{n\rightarrow \infty }P_{n,\theta }\left( \hat{\beta}_{n,i}=0\right) =1%
\text{ \ \ \ holds whenever \ \ \ }\beta _{i}=0\text{.}  \label{part_sparse}
\end{equation}%
E.g., $\hat{\theta}_{n}$ could be a post-model-selection estimator based on
a consistent model selection procedure that only subjects the elements in $%
\beta $ to selection, the elements in $\alpha $ being `protected'.

If we are now interested in a confidence set for $\beta $ that is based on $%
\hat{\beta}_{n}$, we can immediately apply the results obtained sofar: By
viewing $\alpha $ as a `nuisance' parameter, we can use Remark \ref{nuis} to
conclude that Theorem \ref{Th1} applies mutatis mutandis to this situation.
Moreover, combining the reasoning in Remarks \ref{nuis} and \ref{funct}, we
can then immediately obtain a result similar to (\ref{lin_funct}) for
confidence sets for $A\beta $ that are based on $A\hat{\beta}_{n}$, $A$
being an arbitrary matrix of full row rank. For the sake of brevity we do
not spell out the details which are easily obtained from the outline just
given.

The above results, however, do not cover the case where one is interested in
a confidence set for $\theta $ based on a partially sparse estimator $\hat{%
\theta}_{n}$, or more generally the case of confidence sets for $A\theta $
based on $A\hat{\theta}_{n}$, where the linear function $A\theta $ is also
allowed to depend on $\alpha $. For this case we have the following result.

\begin{theorem}
\label{Th2}Suppose the statistical experiment given in (\ref{0}) is such
that for some $\alpha \in \mathbb{R}^{k-k_{\beta }}$ the sequence $%
P_{n,(\alpha ^{\prime },\gamma ^{\prime }/\sqrt{n})^{\prime }}$ is
contiguous w.r.t. $P_{n,(\alpha ^{\prime },0)^{\prime }}$ for every $\gamma
\in \mathbb{R}^{k_{\beta }}$. Let $\hat{\theta}_{n}$ be an estimator
sequence that is partially sparse in the sense of (\ref{part_sparse}). Let $%
A $ be a $q\times k$ matrix of full row rank, which is partitioned
conformably with $\theta $ as $A=(A_{1},A_{2})$, and that satisfies $%
\limfunc{rank}A_{1}<q$. Let $C_{n}$ be a sequence of random sets \textit{%
based on }$A\hat{\theta}_{n}$ (in the sense that $P_{n,\theta }\left( A\hat{%
\theta}_{n}\in C_{n}\right) =1$ for all $\theta \in \mathbb{R}^{k}$). Assume
that $C_{n}$ is a confidence set for $A\theta $ with asymptotic infimal
coverage probability $\delta $, i.e.,%
\[
\delta =\liminf_{n\rightarrow \infty }\inf_{\theta \in \mathbb{R}%
^{k}}P_{n,\theta }\left( A\theta \in C_{n}\right) . 
\]%
Then for every $t\geq 0$ we have%
\begin{equation}
\liminf_{n\rightarrow \infty }\sup_{\theta \in \mathbb{R}^{k}}P_{n,\theta
}\left( \sqrt{n}\limfunc{diam}(C_{n})\geq t\right) \geq \delta .
\label{diam_2}
\end{equation}%
[If the set inside of the probability in (\ref{diam_2}) is not measurable,
the probability is to be replaced by inner probability.]
\end{theorem}

\begin{proof}
Consider sequences $\theta _{n}=(\alpha ^{\prime },\gamma ^{\prime }/\sqrt{n}%
)^{\prime }\in \mathbb{R}^{k}$ where $\alpha $ is as in the theorem. Then
similar as in the proof of Theorem \ref{Th1} exploiting partial sparsity and
contiguity we arrive at%
\begin{eqnarray}
\delta &\leq &\liminf_{n\rightarrow \infty }P_{n,\theta _{n}}\left( A\theta
_{n}\in C_{n}\right)  \nonumber \\
&\leq &\liminf_{n\rightarrow \infty }P_{n,\theta _{n}}\left( A\theta _{n}\in
C_{n},A\hat{\theta}_{n}\in C_{n},\hat{\beta}_{n}=0\right)  \nonumber \\
&\leq &\liminf_{n\rightarrow \infty }P_{n,\theta _{n}}\left( \limfunc{diam}%
(C_{n})\geq \left\Vert A((\alpha -\hat{\alpha}_{n})^{\prime },\gamma
^{\prime }/\sqrt{n})^{\prime }\right\Vert \right) .  \label{ineq}
\end{eqnarray}%
By the assumption on $A$ there exists a vector $\gamma _{0}$ such that $%
A_{2}\gamma _{0}$ is non-zero and is linearly independent of the range space
of $A_{1}$. Consequently, $\Pi A_{2}\gamma _{0}\neq 0$, where $\Pi $ denotes
the orthogonal projection on the orthogonal complement of the range space of 
$A_{1}$. Set $\gamma =c\gamma _{0}$ for arbitrary $c$. Then 
\begin{eqnarray*}
\left\Vert A((\alpha -\hat{\alpha}_{n})^{\prime },\gamma ^{\prime }/\sqrt{n}%
)^{\prime }\right\Vert ^{2} &=&\left\Vert A_{1}(\alpha -\hat{\alpha}%
_{n})+A_{2}\gamma /\sqrt{n}\right\Vert ^{2} \\
&\geq &n^{-1}c^{2}\left\Vert \Pi A_{2}\gamma _{0}\right\Vert ^{2}.
\end{eqnarray*}%
Combined with (\ref{ineq}), this gives 
\[
\delta \leq \liminf_{n\rightarrow \infty }P_{n,\theta _{n}}\left( \sqrt{n}%
\limfunc{diam}(C_{n})\geq \left\vert c\right\vert \left\Vert \Pi A_{2}\gamma
_{0}\right\Vert \right) . 
\]%
Since $\left\Vert \Pi A_{2}\gamma _{0}\right\Vert >0$ by construction and
since $c$ was arbitrary, the result (\ref{diam_2})\ follows upon identifying 
$t$ and $\left\vert c\right\vert \left\Vert \Pi A_{2}\gamma _{0}\right\Vert $%
.
\end{proof}

Some simple generalizations are possible: Inspection of the proof shows that 
$\delta $ may be replaced by $\delta (\alpha )=\liminf_{n\rightarrow \infty
}\inf_{\beta \in \mathbb{R}^{k}}P_{n,\theta }\left( A\theta \in C_{n}\right) 
$ where $\alpha $ is as in the theorem and $\theta =(\alpha ^{\prime },\beta
^{\prime })^{\prime }$. Furthermore, the partial sparsity condition (\ref%
{part_sparse}) only needs to hold for all $\theta =(\alpha ^{\prime },\beta
^{\prime })^{\prime }$ with $\alpha $ as in the theorem. A similar remark
applies to Theorem \ref{Th3} given below.

The condition on $A$ in the above theorem is, for example, satisfied when
considering confidence sets for the entire vector $\theta $ as this
corresponds to the case $A=I_{k}$ (and $q=k$). [The condition is also
satisfied in case $A=(0_{k_{\beta }\times (k-k_{\beta })},I_{k_{\beta }})$
which corresponds to setting confidence sets for $\beta $. However, in this
case already the extension of Theorem \ref{Th1} discussed prior to Theorem %
\ref{Th2} applies.]

Theorem \ref{Th2} does not cover the case where a confidence set is desired
for $\alpha $ only (i.e., $A=(I_{k-k_{\beta }},0_{(k-k_{\beta })\times
k_{\beta }})$). In fact, without further assumptions on the estimator $\hat{%
\theta}_{n}$ no result of the above sort is in general possible in this case
(to see this consider the case where $\hat{\alpha}_{n}$ and $\hat{\beta}_{n}$
are independent and $\hat{\alpha}_{n}$ is a well-behaved estimator).
However, under additional assumptions, results that show that confidence
sets for $\alpha $ are also necessarily large will be obtained next. We
first present the result and subsequently discuss the assumptions.

\begin{theorem}
\label{Th3}Suppose the statistical experiment given in (\ref{0}) is such
that for some $\alpha \in \mathbb{R}^{k-k_{\beta }}$ the sequence $%
P_{n,(\alpha ^{\prime },\gamma ^{\prime }/\sqrt{n})^{\prime }}$ is
contiguous w.r.t. $P_{n,(\alpha ^{\prime },0)^{\prime }}$ for every $\gamma
\in \mathbb{R}^{k_{\beta }}$. Let $\hat{\theta}_{n}$ be an estimator
sequence that is partially sparse in the sense of (\ref{part_sparse}).
Suppose that there exists a $(k-k_{\beta })\times k_{\beta }$-matrix $D$
such that for every $\gamma $ the random vector $n^{1/2}(\hat{\alpha}%
_{n}-\alpha )$ converges in $P_{n,(\alpha ^{\prime },\gamma ^{\prime }/\sqrt{%
n})^{\prime }}$-distribution to $Z+D\gamma $ where $Z$ is a $(k-k_{\beta
})\times 1$\ random vector with a distribution that is independent of $%
\gamma $. Let $A$ be a $q\times k$ matrix of full row rank, which is
partitioned conformably with $\theta $ as $A=(A_{1},A_{2})$, and assume that 
$A_{1}D-A_{2}\neq 0$. Let $C_{n}$ be a sequence of random sets \textit{based
on }$A\hat{\theta}_{n}$ (in the sense that $P_{n,\theta }\left( A\hat{\theta}%
_{n}\in C_{n}\right) =1$ for all $\theta \in \mathbb{R}^{k}$). Assume that $%
C_{n}$ is a confidence set for $A\theta $ with asymptotic infimal coverage
probability $\delta $, i.e.,%
\[
\delta =\liminf_{n\rightarrow \infty }\inf_{\theta \in \mathbb{R}%
^{k}}P_{n,\theta }\left( A\theta \in C_{n}\right) . 
\]%
Then for every $t\geq 0$ we have%
\begin{equation}
\liminf_{n\rightarrow \infty }\sup_{\theta \in \mathbb{R}^{k}}P_{n,\theta
}\left( \sqrt{n}\limfunc{diam}(C_{n})\geq t\right) \geq \delta .
\label{diam_3}
\end{equation}%
[If the set inside of the probability in (\ref{diam_3}) is not measurable,
the probability is to be replaced by inner probability.]
\end{theorem}

\begin{proof}
Consider sequences $\theta _{n}=(\alpha ^{\prime },\gamma ^{\prime }/\sqrt{n}%
)^{\prime }\in \mathbb{R}^{k}$ where $\alpha $ is as in the theorem. Then
for every $t\geq 0$ we have%
\begin{eqnarray}
\delta &\leq &\liminf_{n\rightarrow \infty }P_{n,\theta _{n}}\left( A\theta
_{n}\in C_{n}\right) =\liminf_{n\rightarrow \infty }P_{n,\theta _{n}}\left(
A\theta _{n}\in C_{n},A\hat{\theta}_{n}\in C_{n}\right)  \nonumber \\
&\leq &\liminf_{n\rightarrow \infty }P_{n,\theta _{n}}\left( A\theta _{n}\in
C_{n},A\hat{\theta}_{n}\in C_{n},n^{1/2}\left\Vert A(\hat{\theta}_{n}-\theta
_{n})\right\Vert \geq t\right)  \nonumber \\
&&+\limsup_{n\rightarrow \infty }P_{n,\theta _{n}}\left( n^{1/2}\left\Vert A(%
\hat{\theta}_{n}-\theta _{n})\right\Vert <t\right)  \nonumber \\
&\leq &\liminf_{n\rightarrow \infty }P_{n,\theta _{n}}\left( n^{1/2}\limfunc{%
diam}(C_{n})\geq t\right)  \nonumber \\
&&+\limsup_{n\rightarrow \infty }P_{n,\theta _{n}}\left( n^{1/2}\left\Vert A(%
\hat{\theta}_{n}-\theta _{n})\right\Vert <t\right) .  \label{cruc}
\end{eqnarray}%
Exploiting partial sparsity and contiguity we get%
\begin{eqnarray}
&&\limsup_{n\rightarrow \infty }P_{n,\theta _{n}}\left( n^{1/2}\left\Vert A(%
\hat{\theta}_{n}-\theta _{n})\right\Vert <t\right)  \nonumber \\
&\leq &\limsup_{n\rightarrow \infty }P_{n,\theta _{n}}\left( \hat{\beta}%
_{n}=0,n^{1/2}\left\Vert A(\hat{\theta}_{n}-\theta _{n})\right\Vert <t\right)
\nonumber \\
&&+\limsup_{n\rightarrow \infty }P_{n,\theta _{n}}\left( \hat{\beta}_{n}\neq
0\right)  \nonumber \\
&=&\limsup_{n\rightarrow \infty }P_{n,\theta _{n}}\left( \hat{\beta}%
_{n}=0,n^{1/2}\left\Vert A(\hat{\theta}_{n}-\theta _{n})\right\Vert <t\right)
\nonumber \\
&\leq &\limsup_{n\rightarrow \infty }P_{n,\theta _{n}}\left(
n^{1/2}\left\Vert A_{1}(\hat{\alpha}_{n}-\alpha )-A_{2}\gamma /\sqrt{n}%
\right\Vert <t\right)  \nonumber \\
&=&\limsup_{n\rightarrow \infty }P_{n,\theta _{n}}\left( \left\Vert
X_{n}+(A_{1}D-A_{2})\gamma \right\Vert <t\right)  \nonumber \\
&\leq &\limsup_{n\rightarrow \infty }P_{n,\theta _{n}}\left( \left\Vert
X_{n}\right\Vert >\left\Vert (A_{1}D-A_{2})\gamma \right\Vert -t\right)
\label{ine}
\end{eqnarray}%
where $X_{n}$ converges to $A_{1}Z$ in $P_{n,\theta _{n}}$-distribution.
Since $A_{1}D-A_{2}\neq 0$ by assumption, we can find a $\gamma $ such that $%
\left\Vert (A_{1}D-A_{2})\gamma \right\Vert -t$ is arbitrarily large, making
the far right-hand side of (\ref{ine}) arbitrarily small. This, together
with (\ref{cruc}), establishes the result.
\end{proof}

Note that the case where a confidence set for $\alpha $ is sought, that is, $%
A=(I_{k-k_{\beta }},0_{(k-k_{\beta })\times k_{\beta }})$, which was not
covered by Theorem \ref{Th2}, is covered by Theorem \ref{Th3} except in the
special case where $D=0$.

The weak convergence assumption in the above theorem merits some discussion:
Suppose $\hat{\theta}_{n}$ is a post-model-selection estimator based on a
model selection procedure that consistently finds the zeroes in $\beta $ and
then computes $\hat{\theta}_{n}$ as the restricted maximum likelihood
estimator $\hat{\theta}_{n}(R)$ under the zero-restrictions in $\beta $.
Under the usual regularity conditions, the restricted maximum likelihood
estimator $\hat{\alpha}_{n}(R)$ for $\alpha $ will then satisfy that $%
n^{1/2}(\hat{\alpha}_{n}(R)-\alpha )$ converges to a $N(D\gamma ,\Sigma )$%
-distribution under the sequence of local alternatives $\theta _{n}=(\alpha
^{\prime },\gamma ^{\prime }/\sqrt{n})^{\prime }$. Since $\lim_{n\rightarrow
\infty }P_{n,\theta _{n}}\left( \hat{\beta}_{n}=0\right) =1$ by partial
sparsity and contiguity, the estimators $\hat{\alpha}_{n}$ and $\hat{\alpha}%
_{n}(R)$ coincide with $P_{n,\theta _{n}}$-probability approaching one. This
shows that the assumption on $\hat{\alpha}_{n}$ will typically be satisfied
for such post-model-selection estimators with $Z\sim N(0,\Sigma )$. [For a
precise statement of such a result in a simple example see Leeb and P\"{o}%
tscher (2005, Proposition A.2).] While we expect that this assumption on the
asymptotic behavior of $\hat{\alpha}_{n}$ is also shared by many other
partially sparse estimators, this remains to be verified on a case by case
basis.

\section{An Example: A confidence set based on a hard-thresholding estimator}

Suppose the data $y_{1},\ldots ,y_{n}$ are independent identically
distributed as $N(\theta ,1)$, $\theta \in \mathbb{R}$. Let the
hard-thresholding estimator $\hat{\theta}_{n}$ be given by%
\[
\hat{\theta}_{n}=\bar{y}\boldsymbol{1}(\left\vert \bar{y}\right\vert >\eta
_{n}) 
\]%
where the threshold $\eta _{n}$ is a positive real number and $\bar{y}$
denotes the maximum likelihood estimator, i.e., the arithmetic mean of the
data. Of course, $\hat{\theta}_{n}$ is nothing else than a
post-model-selection estimator following a $t$-type test of the hypothesis $%
\theta =0$ versus the alternative $\theta \neq 0$. It is well-known and easy
to see that $\hat{\theta}_{n}$ satisfies the sparsity condition if $\eta
_{n}\rightarrow 0$ and $n^{1/2}\eta _{n}\rightarrow \infty $ (i.e., the
underlying model selection procedure is consistent); in this case then $%
n^{1/2}(\hat{\theta}_{n}-\theta )$ converges to a standard normal
distribution if $\theta \neq 0$, whereas it converges to pointmass at zero
if $\theta =0$. Note that $\hat{\theta}_{n}$ -- with such a choice of the
threshold $\eta _{n}$ -- is an instance of Hodges' estimator. In contrast,
if $\eta _{n}\rightarrow 0$ and $n^{1/2}\eta _{n}\rightarrow e$, $0\leq
e<\infty $, the estimator $\hat{\theta}_{n}$ is a post-model-selection
estimator based on a conservative model selection procedure. See P\"{o}%
tscher and Leeb (2007) for further discussion and references.

In the consistent model selection case the estimator possesses the
\textquotedblleft oracle\textquotedblright\ property suggesting as a
confidence interval the \textquotedblleft naive\textquotedblright\ interval
given by $C_{n}^{naive}=\{0\}$ if $\hat{\theta}_{n}=0$ and by $%
C_{n}^{naive}=[\hat{\theta}_{n}-z_{(1-\delta )/2},\hat{\theta}%
_{n}+z_{(1-\delta )/2}]$ otherwise, where $\delta $ is the nominal coverage
level and $z_{(1-\delta )/2}$ is the $1-(1-\delta )/2$-quantile of the
standard normal distribution. This interval satisfies $P_{n,\theta }(\theta
\in C_{n}^{naive})\rightarrow \delta $ for every $\theta $, but -- as
discussed in the introduction and as follows from the results in Section 2
-- it is not honest and, in fact, has infimal coverage probability
converging to zero. A related, but infeasible, construction is to consider
the intervals $C_{n}^{\ast }=[\hat{\theta}_{n}-c_{n}(\theta ),\hat{\theta}%
_{n}+c_{n}(\theta )]$ where $c_{n}(\theta )$ is chosen as small as possible
subject to $P_{n,\theta }(\theta \in C_{n}^{\ast })=\delta $ for every $%
\theta $. [Note that $C_{n}^{naive}$ can be viewed as being obtained from $%
C_{n}^{\ast }$ by replacing $c_{n}(\theta )$ by the limits $c_{\infty
}(\theta )$ for $n\rightarrow \infty $, where $c_{\infty }(\theta )=0$ if $%
\theta =0$ and $c_{\infty }(\theta )=z_{(1-\delta )/2}$ if $\theta \neq 0$,
and then by replacing $\theta $ by $\hat{\theta}_{n}$ in $c_{\infty }(\theta
)$.] An obvious idea to obtain a feasible and honest interval is now to use $%
c_{n}=\max_{\theta \in \mathbb{R}}$ $c_{n}(\theta )$ as the half-length of
the interval, i.e. $C_{n}=[\hat{\theta}_{n}-c_{n},\hat{\theta}_{n}+c_{n}]$.
From Theorem \ref{Th1} we know that $\sqrt{n}c_{n}\rightarrow \infty $ in
the case where $\eta _{n}\rightarrow 0$ and $n^{1/2}\eta _{n}\rightarrow
\infty $ (and if $\delta >0$), but it is instructive to study the behavior
of $C_{n}$ in more detail.

We therefore consider now confidence intervals $C_{n}$ for $\theta $ of the
form $C_{n}=[\hat{\theta}_{n}-a_{n},\hat{\theta}_{n}+b_{n}]$ with
nonnegative constants $a_{n}$ and $b_{n}$ (thus removing the symmetry
restriction on the interval). Note that the subsequent result is a
finite-sample result and hence does not involve any assumptions on the
behavior of $\eta _{n}$.

\begin{proposition}
For every $n\geq 1$, the interval $C_{n}=[\hat{\theta}_{n}-a_{n},\hat{\theta}%
_{n}+b_{n}]$ has an infimal coverage probability satisfying 
\begin{eqnarray*}
&&\inf_{\theta \in \mathbb{R}}P_{n,\theta }\left( \theta \in C_{n}\right)  \\
&=&\left\{ 
\begin{array}{ll}
\Phi (n^{1/2}(a_{n}-\eta _{n}))-\Phi (-n^{1/2}b_{n}) & \text{if \ \ }\eta
_{n}\leq a_{n}+b_{n}\text{ \ and \ }a_{n}\leq b_{n} \\ 
\Phi (n^{1/2}a_{n})-\Phi (n^{1/2}(-b_{n}+\eta _{n})) & \text{if \ \ }\eta
_{n}\leq a_{n}+b_{n}\text{ \ and \ }a_{n}\geq b_{n} \\ 
0 & \text{if \ \ \ }\eta _{n}>a_{n}+b_{n}%
\end{array}%
\right. ,
\end{eqnarray*}%
where $\Phi $ denotes the standard normal cumulative distribution function.
\end{proposition}

\begin{proof}
Elementary calculations and the fact that $n^{1/2}(\bar{y}-\theta )$ is
standard normally distributed give for the coverage probability $%
p_{n}(\theta )=P_{n,\theta }\left( \theta \in C_{n}\right) $%
\begin{eqnarray*}
p_{n}(\theta ) &=&P_{n,\theta }\left( -n^{1/2}b_{n}\leq n^{1/2}(\hat{\theta}%
_{n}-\theta )\leq n^{1/2}a_{n}\right) \\
&=&\Pr \left( -n^{1/2}b_{n}\leq Z\leq n^{1/2}a_{n},\left\vert
Z+n^{1/2}\theta \right\vert >n^{1/2}\eta _{n}\right) \\
&&+\Pr \left( -b_{n}\leq -\theta \leq a_{n},\left\vert Z+n^{1/2}\theta
\right\vert \leq n^{1/2}\eta _{n}\right) ,
\end{eqnarray*}%
where $Z$ is a standard normally distributed random variable and $\Pr $
denotes a generic probability. Simple, albeit tedious computations give the
coverage probability as follows. If $\eta _{n}>a_{n}+b_{n}$%
\[
p_{n}(\theta )=\left\{ 
\begin{array}{ll}
\Phi (n^{1/2}a_{n})-\Phi (-n^{1/2}b_{n}) & \text{if \ }\theta <-a_{n}-\eta
_{n}\text{ \ or \ }\theta >b_{n}+\eta _{n} \\ 
\Phi (n^{1/2}(-\theta -\eta _{n}))-\Phi (-n^{1/2}b_{n}) & \text{if \ }%
-a_{n}-\eta _{n}\leq \theta <b_{n}-\eta _{n} \\ 
0 & \text{if\ \ }b_{n}-\eta _{n}\leq \theta <-a_{n}\text{ \ or \ }%
b_{n}<\theta \leq -a_{n}+\eta _{n} \\ 
\Phi (n^{1/2}(-\theta +\eta _{n}))-\Phi (n^{1/2}(-\theta -\eta _{n})) & 
\text{if\ \ }-a_{n}\leq \theta \leq b_{n} \\ 
\Phi (n^{1/2}a_{n})-\Phi (n^{1/2}(-\theta +\eta _{n})) & \text{if\ \ }%
-a_{n}+\eta _{n}<\theta \leq b_{n}+\eta _{n}%
\end{array}%
\right. . 
\]%
Hence, the infimal coverage probability in this case is obviously zero.
Next, if $(a_{n}+b_{n})/2\leq \eta _{n}\leq a_{n}+b_{n}$ then%
\[
p_{n}(\theta )=\left\{ 
\begin{array}{ll}
\Phi (n^{1/2}a_{n})-\Phi (-n^{1/2}b_{n}) & \text{if \ }\theta <-a_{n}-\eta
_{n}\text{ \ or \ }\theta >b_{n}+\eta _{n} \\ 
\Phi (n^{1/2}(-\theta -\eta _{n}))-\Phi (-n^{1/2}b_{n}) & \text{if \ }%
-a_{n}-\eta _{n}\leq \theta <-a_{n} \\ 
\Phi (n^{1/2}(-\theta +\eta _{n}))-\Phi (-n^{1/2}b_{n}) & \text{if\ \ }%
-a_{n}\leq \theta <b_{n}-\eta _{n} \\ 
\Phi (n^{1/2}(-\theta +\eta _{n}))-\Phi (n^{1/2}(-\theta -\eta _{n})) & 
\text{if\ \ }b_{n}-\eta _{n}\leq \theta \leq -a_{n}+\eta _{n} \\ 
\Phi (n^{1/2}a_{n})-\Phi (n^{1/2}(-\theta -\eta _{n})) & \text{if\ \ }%
-a_{n}+\eta _{n}<\theta \leq b_{n} \\ 
\Phi (n^{1/2}a_{n})-\Phi (n^{1/2}(-\theta +\eta _{n})) & \text{if\ \ }%
b_{n}<\theta \leq b_{n}+\eta _{n}%
\end{array}%
\right. , 
\]%
and if $\eta _{n}<(a_{n}+b_{n})/2$%
\[
p_{n}(\theta )=\left\{ 
\begin{array}{ll}
\Phi (n^{1/2}a_{n})-\Phi (-n^{1/2}b_{n}) & 
\begin{array}{l}
\text{if \ }\theta <-a_{n}-\eta _{n}\text{ \ or \ }\theta >b_{n}+\eta _{n}
\\ 
\text{or \ }-a_{n}+\eta _{n}\leq \theta \leq b_{n}-\eta _{n}%
\end{array}
\\ 
\Phi (n^{1/2}(-\theta -\eta _{n}))-\Phi (-n^{1/2}b_{n}) & \text{ if \ }%
-a_{n}-\eta _{n}\leq \theta <-a_{n} \\ 
\Phi (n^{1/2}(-\theta +\eta _{n}))-\Phi (-n^{1/2}b_{n}) & \text{ if \ }%
-a_{n}\leq \theta <-a_{n}+\eta _{n} \\ 
\Phi (n^{1/2}a_{n})-\Phi (n^{1/2}(-\theta -\eta _{n})) & \text{ if\ \ }%
b_{n}-\eta _{n}<\theta \leq b_{n} \\ 
\Phi (n^{1/2}a_{n})-\Phi (n^{1/2}(-\theta +\eta _{n})) & \text{ if\ \ }%
b_{n}<\theta \leq b_{n}+\eta _{n}%
\end{array}%
\right. . 
\]%
Inspection shows that in both cases the function does not have a minimum,
but the infimum equals the smaller of the left-hand side limit $%
p_{n}(-a_{n}-)$ and the right-hand side limit $p_{n}(b_{n}+)$, which shows
that the infimum of $p_{n}(\theta )$ equals $\min [\Phi (n^{1/2}(a_{n}-\eta
_{n}))-\Phi (-n^{1/2}b_{n}),\Phi (n^{1/2}a_{n})-\Phi (n^{1/2}(-b_{n}+\eta
_{n}))]$.
\end{proof}

As a point of interest we note that the coverage probability $p_{n}(\theta )$
has exactly two discontinuity points (jumps), one at $\theta =-a_{n}$ and
one at $\theta =b_{n}$, except in the trivial case $a_{n}=b_{n}=0$ where the
two discontinuity points merge into one.

An immediate consequence of the above proposition is that $n^{1/2}\limfunc{%
diam}(C_{n})=n^{1/2}(a_{n}+b_{n})$ is not less than $n^{1/2}\eta _{n}$,
provided the infimal coverage probability is positive. Hence, in case that $%
\eta _{n}\rightarrow 0$ and $n^{1/2}\eta _{n}\rightarrow \infty $, i.e., in
case that $\hat{\theta}_{n}$ is sparse, we see that $n^{1/2}\limfunc{diam}%
(C_{n})\rightarrow \infty $, which of course just confirms the general
result obtained in Theorem \ref{Th1} above. [In fact, this result is a bit
stronger as only the infimal coverage probabilities need to be positive, and
not their limes inferior.]

If the interval is symmetric, i.e., $a_{n}=b_{n}$ holds, and $a_{n}\geq \eta
_{n}/2$ is satisfied, the infimal coverage probability becomes $\Phi
(n^{1/2}a_{n})-\Phi (n^{1/2}(-a_{n}+\eta _{n}))$. Since this expression is
zero if $a_{n}=\eta _{n}/2$, and is strictly increasing to one as $a_{n}$
goes to infinity, any prescribed infimal coverage probability less than one
is attainable. Suppose $0<\delta <1$ is given. Then the (shortest)
confidence interval $C_{n}$ of the form $[\hat{\theta}_{n}-a_{n},\hat{\theta}%
_{n}+a_{n}]$ with infimal coverage probability equal to $\delta $ has to
satisfy $a_{n}\geq \eta _{n}/2$ and%
\[
\Phi (n^{1/2}a_{n})-\Phi (n^{1/2}(-a_{n}+\eta _{n}))=\delta . 
\]%
If now $\eta _{n}\rightarrow 0$ and $n^{1/2}\eta _{n}\rightarrow \infty $,
i.e., if $\hat{\theta}_{n}$ is sparse, it follows that $n^{1/2}a_{n}%
\rightarrow \infty $ and%
\[
n^{1/2}(-a_{n}+\eta _{n})\rightarrow \Phi ^{-1}(1-\delta ) 
\]%
or in other words that $a_{n}\geq \eta _{n}/2$ has to satisfy%
\begin{equation}
a_{n}=\eta _{n}-n^{-1/2}\Phi ^{-1}(1-\delta )+o(n^{-1/2}).  \label{A}
\end{equation}%
Conversely, any $a_{n}\geq \eta _{n}/2$ satisfying (\ref{A}) generates a
confidence interval with asymptotic infimal coverage probability equal to $%
\delta $. We observe that (\ref{A}) shows that $\kappa _{n}\limfunc{diam}%
(C_{n})=2\kappa _{n}a_{n}\rightarrow \infty $ for any sequence that
satisfies $\kappa _{n}\eta _{n}\rightarrow \infty $, which includes
sequences that are $o(n^{1/2})$ by the assumptions on $\eta _{n}$. Hence,
this result is stronger than what is obtained from applying Theorem \ref{Th1}
(or its Corollary) to this example, and illustrates the discussion in Remark %
\ref{rate}.

\end{document}